\documentclass[11pt]{article}
\usepackage[margin=1in]{geometry}
\usepackage{amsmath,amssymb,amsthm,mathtools}
\usepackage{enumitem}
\usepackage[colorlinks=true,linkcolor=blue,citecolor=blue,urlcolor=blue]{hyperref}
\usepackage{thmtools}
\usepackage{thm-restate}
\usepackage{authblk}

\newtheorem{theorem}{Theorem}[section]
\newtheorem{definition}[theorem]{Definition}
\newtheorem{lemma}[theorem]{Lemma}

\newtheorem{observation}[theorem]{Observation}

\DeclareMathOperator{\A}{\mathcal{A}}
\DeclareMathOperator{\boxi}{box}

\DeclareMathOperator{\dimt}{dim_{TH}}
\DeclareMathOperator{\dimc}{dim_{C4F}}
\DeclareMathOperator{\dimi}{dim_{INT}}
\DeclareMathOperator{\dima}{dim_{\mathcal{A}}}
\DeclareMathOperator{\cova}{cov_{\mathcal{A}}}
\DeclareMathOperator{\covk}{cov_{2K2F}}
\DeclareMathOperator{\covt}{cov_{TH}}

\newcommand{\bbR}{\mathbb R}
\newcommand{\calL}{\mathcal L}
\def\G#1{G=G_1\cap G_2\cap\cdots\cap G_{#1}}
\title{New bounds on the boxicity of Kneser graphs}
\author[1]{L. Sunil Chandran%
  \thanks{Email: \texttt{sunil@iisc.ac.in}}}

\author[1]{Suraj Kumar Sahoo%
  \thanks{Email: \texttt{surajks@iisc.ac.in}}}
\affil[1]{Indian Institute of Science}
\date{}
\begin{document}
\maketitle
\begin{abstract} 
For a graph $G=(V,E)$, the boxicity $\operatorname{box}(G)$ is the least integer $d$ such that
 $G$ is the intersection graph of axis-parallel boxes in $\mathbb{R}^d$. The Kneser graph $K(n,k)$ is the graph whose vertex set corresponds to the family of $k$ sized subsets of an $n$ sized ground set, and a pair of vertices are adjacent if and only if the corresponding sets are disjoint. 

In this short paper, we provide new bounds for the boxicity of $K(n,k)$. 
We first show that for $k\ge 2$ and $n\ge 2k+1$, $\operatorname{box}(K(n,k))\ge n-2k$, thereby improving the lower bounds 
of Caoduro and Lichev [Discrete Mathematics 346(5):113333 (2023)] when $n\geq 3k+1$ and $k\ge 3$. The proof of this result uses the topological method and follows a strategy similar to Greene's proof of the Kneser's conjecture on the chromatic number of $K(n,k)$.

We also show that for $k\ge 2$, $t\ge 1$ and $n=2k+t$, $\operatorname{box}(K(2k+t,k)) \leq  O( t^3\log(2k+t) ).$
In particular, for every $t=t(k)\ge1$ with $t^3\log k=o(k)$, $\operatorname{box}(K(2k+t,k))=o(k)$. This is the first upper bound on the boxicity of $K(n,k)$ that is sublinear in $n$. We also show that for large enough $k$, $\operatorname{box}(K(2k+1,k))\le \frac{3\ln (2k+1)}{\ln\ln (2k+1)}$.
\end{abstract}

\section{Introduction}
In this paper all graphs are finite, simple and undirected. We will denote a graph $G$ defined on the vertex set $V$ and edge set $E$ as $G=(V,E).$ 

For a graph $G=(V,E)$, the \emph{boxicity} $\boxi(G)$ is the least integer $d$ such that
$G$ is the intersection graph of axis-parallel boxes in $\bbR^d$. If $G$ is a complete graph, we take $\boxi(G)=0$ by convention. 

For a graph $G$, several bounds on $\boxi(G)$ have been established, notably in terms of its treewidth, vertex cover number, maximum degree, and degeneracy of $G$~\cite{chandran2007boxicity, chandran2009cubicity,ADIGA20142}. Bounds have also been derived when $G$ belongs to special classes of graphs like planar graphs, chordal graphs, AT-free graphs, and co-comparability graphs~\cite{thomassen1986interval, bhowmick2010boxicity,leaf}. For more on boxicity of graphs, readers are referred to the recent survey~\cite{chandran2026survey}.

The Kneser graph $K(n,k)$ is the graph whose vertex set corresponds to the family of $k$ sized subsets of an $n$ sized ground set, and a pair of vertices are adjacent if and only if the corresponding sets are disjoint. The study of the boxicity of the Kneser graph $K(n,k)$ was initiated by Caoduro and Lichev~\cite{CAODURO2023113333} who showed that for $n\ge 2k+1$, $\boxi(K(n,k))\leq n-2$ and for $n\geq 2k^3-2k^2+1$, $\boxi(K(n,k))\geq n-\frac{1}{2}(13k^2-11k+16)$. They also showed that $\boxi(K(n,k))\ge \frac{1}{2}(n-2k-2)$ for $2k+1\le n\le 3k$ and $\boxi(K(n,k))\ge\frac{1}{2}(n-k-5)$ for $ (n\ge 3k+1).$ Later Caoduro and Seb{\H{o}}~\cite{caoduro2023boxicityy} showed that for $n\ge 5$, $\boxi(K(n,2))=n-2$ and conjectured whether for $n\geq 2k+1$, $\boxi(K(n,k))=n-k$. This conjecture was refuted in~\cite{chandran2026geometric} where the exact boxicity of $K(7,3)$ was computationally shown to be 3.   

The lower bound $\boxi(K(n,k))\geq n-\frac{1}{2}(13k^2-11k+16)$ for $n\geq 2k^3-2k^2+1$ established in~\cite{CAODURO2023113333} is achieved by using ideas based on the vertex boundaries of a graph as developed in~\cite{adiga2014lower}. Despite its promise, a drawback of this method is that it often involves tedious computations of the vertex boundaries for each subset of vertices. Furthermore, it is a natural question to ask whether or not existing methods can be used to get better bounds for $\boxi(K(n,k))$. We answer this question in the positive. 
\subsection{Our results}
For ease of presentation, we state below the definition of`the intersection dimension' of a graph, a notion introduced in~\cite{Kratochvil}. For $\ell\geq 0$ and a collection of graphs $\{G_i=(V,E_i)\}_{1\le i\le\ell}$, the graph $\bigcap_{1\le i\le\ell}G_i$ is the graph defined on $V$ as its set of vertices such that its edge set is $\bigcap_{1\le i\le\ell}E_i$.
\begin{definition}[\cite{Kratochvil}]\label{defdim}
    Let $G=(V,E)$ be a graph and $\mathcal{A}$ be a class of graphs that contains the complete graph $K_n$ and $K_n\setminus e$, for an edge e, $\forall n\in \mathbb{N}$. \emph{The intersection dimension of $G$ with respect to $\mathcal{A}$} (denoted as $\dima(G)$) is defined as follows:
    $$\dima(G)=\min\left\{\ell\colon  \exists\text{ a collection of graphs }\{G_i\in \A\}_{1\le i\le\ell}\text{ such that }G=\bigcap_{i=1}^{\ell}G_i  \right\}.$$
    If $G$ is a complete graph, we take $\dima(G)=0$ by convention.
\end{definition}
Our first result optimally represents $K(n,k)$ as the intersection of $C_4$-free graphs. Let $\operatorname{C4F}$ be the class of $C_4$-free graphs. 
\begin{theorem}\label{c4f}
   For $k\ge 2$ and $n\ge 2k+1$, $\dimc(K(n,k))=n-2k$.
\end{theorem}
This immediately gives us $\dima(K(n,k))\ge n-2k$ when $\A$ is the class of chordal graphs, split graphs, threshold graphs or any other graph class that is $C_4$-free. In particular, this gives us a lower bound for boxicity.
\begin{restatable}{corollary}{thbox}\label{thbox}
    For $k\ge 2$ and $n\ge 2k+1$, $ \boxi(K(n,k))\ge n-2k$.
\end{restatable}
 The proof of Theorem~\ref{c4f} uses a generalization of a theorem by Lusternik, Schnirelman and Borsuk developed by Greene~\cite{Greene01122002} and bears close similarity to his simplified proof of the Kneser's conjecture for the chromatic number of $K(n,k)$ (the first proof of this conjecture was given by Lov\H{a}sz~\cite{LOVASZ1978319} and had been simplified once by B\'ar\'any~\cite{BARANY1978325} before Greene~\cite{Greene01122002}). The strength of this method is illustrated in the clean and unconditional bound obtained, and also in its ease of presentation.  

We also show the first upper bounds for the boxicity of $K(n,k)$ when $n$ is `close' to $2k$, thereby bringing us closer to determining the true behaviour of $\boxi(K(n,k))$ in this regime. We do this by giving upper bounds for the intersection dimension of $K(n,k)$ with respect to threshold graphs. Let $\operatorname{TH}$ be the class of threshold graphs.
 \begin{restatable}{theorem}{sublinear}\label{sublinear}
For $k\ge 2$ and $t\ge 1$, we have $$ \boxi\left(K\left(2k+t,k\right)\right)\le \dimt(K(2k+t,k))\leq  O( t^3\log(2k+t) ).$$
In particular, for every $t=t(k)\ge1$ with $t^3\log k=o(k)$, $\boxi\left(K\left(2k+t,k\right)\right)\le \dimt(K(2k+t,k))=o(k)$.  
\end{restatable}
\noindent When $n=2k+1$, we can do much better. 
 \begin{restatable}{theorem}{odd}\label{odd}\label{odd}
   For $k\ge 2$, $ \boxi(K(2k+1,k))\le\dimt(K(2k+1,k))\le 2\min\{j \colon j! \ge k\} + 2.$ In particular, for large enough $k$, $\boxi(K(2k+1,k))\le \frac{3\ln (2k+1)}{\ln\ln (2k+1)}$. 
\end{restatable}
 Let $B(a,b;c)$ be the poset of height 2 consisting of the subsets of $[c]$ of size $a$ and $b$
ordered by the set-inclusion relation. It was observed in Remark 4.11 of~\cite{CAODURO2023113333} that  $$\boxi(K(n,k))\ge \frac{1}{2}\dim (B(k-1,n-k-1;n-1))-2\geq  \frac{1}{2}\dim (B(1,n-2k+1;n-k+1))-2.$$ As result of Theorem 4.5 in~\cite{kierstead1999dimension}, there is a positive constant $c$ such that $$\boxi(K(2k+t,k))\ge \begin{cases}
 c2^{t+1} \log\log (k+t+1) -2& \text{if }1\le t\le \log\log (k+t+1)-\log\log\log (k+t+1)-1  \\ 
 c(t+1)^2\cdot \frac{\log (t+k+1)}{\log (t+1)}-2& \text{if }1+\log\log (k+t+1)\leq t\le 1+2\sqrt{k+t+1}
\end{cases}. $$
Our upper bounds given in Theorems~\ref{sublinear} and~\ref{odd} complement the aforementioned lower bounds and are, in fact, the only non-trivial upper bounds for boxicity except for the $n-2$ bound given in~\cite{CAODURO2023113333} for the general case.   
\section{Preliminaries}
 We write $\overline G$ for the graph
complement of $G$. For a positive integer $x$, we will use the notation $[x]$ to indicate the set $\{1,2,\ldots, x\}$ with $[0]=\emptyset$. We will denote an interval graph and its interval representation by the same symbol when there is no confusion. For graphs $G$ and $H$, we say that $G$ is $H$-free if $G$ does not contain $H$ as an induced subgraph. 

As alluded to earlier, $\boxi(G)=\dimi(G)$, where INT is the class of interval graphs. We will use this formulation as the working definition for boxicity instead of the standard definition given earlier. 

For ease of understanding of the proofs, we explicitly note down the two criteria that are required to be satisfied for a graph to be an intersection of graphs. For graphs $G=(V,E)$ and $\{G_i=(V,E_i)\}_{i\in [\ell]}$ for some $\ell\ge 0$, we have $\G{\ell}$ if and only if both the conditions below are satisfied.
\begin{enumerate}
    \item $E\subseteq E_i$ for all $i\in [\ell]$, i.e., $G_i$ is a supergraph of $G$ for all $i\in [\ell]$. In other words, if $uv\in E$ then $uv\in E_i$, $\forall i\in [\ell]$.
    \item For distinct vertices $u,v\in V$, if $uv\notin E$ then there exists a $j\in [\ell]$ such that $uv\notin E_j$.
\end{enumerate}

 For graphs $\{G_i=(V,E_i)\}_{i\in [\ell]}$ for some $\ell\ge 0$, we can also define $\bigcup_{i\in [\ell]}G_i$ analogously. There is a notion closely related to the intersection dimension of a graph, known as the \emph{covering number} of the graph. 
\begin{definition}[\cite{Kratochvil}]\label{defcov}
    Let $G=(V,E)$ be a graph and $\mathcal{A}$ be a class of graphs that contains the empty graph $\overline{K_n}$ and $\overline{K_n}\cup e$, for an edge $e$, $\forall n\in \mathbb{N}$. \emph{The covering number of $G$ with respect to $\mathcal{A}$} (denoted as $\cova(G)$) is defined as follows:
    $$\cova(G)=\min\left\{\ell\colon  \exists\text{ a collection of graphs }\{G_i\in \A\}_{i\in [\ell]}\text{ such that }G=\bigcup_{i=1}^{\ell}G_i  \right\}.$$
    If $G$ is an empty graph, we take $\cova(G)=0$ by convention.
\end{definition}
It is clear that for a graph $G$ and graph class $\A$, $\cova(G)=\operatorname{\dim_{co-\A}}(\overline{G})$, where co-$\A=\{H\colon \overline{H}\in \A\}$. 

\noindent The following observation is useful and we state it without a proof.
\begin{observation}\label{sub}
Let $\A$ be closed under vertex deletion. Then for a graph $G$ and an induced subgraph $H$ of $G$, $\dima(G)\geq \dima(H)$.
\end{observation}
Recall that we denote the class of threshold graphs as $\operatorname{TH}$. We formally define this class of graphs. 
\begin{definition}
    Let $s\in \mathbb{R}_{\geq 0}$ and let $V$ be a finite set. Let $w:V\rightarrow \mathbb{R}_{\geq 0}$ be a function. We define $\Gamma(s,w)$ to be the graph such that the vertex set of $\Gamma(s,w)$ is $V$ and for distinct vertices $u,v\in V$, $uv$ is an edge in $\Gamma(s,w)$ if and only if $w(u) + w(v) \geq  s$. A graph $G=(V,E)$ is said to be a threshold graph if for some $s\in \mathbb{R}_{\geq 0}$ and $w:V\rightarrow \mathbb{R}_{\geq 0}$, $G$ is isomorphic to $\Gamma(s,w)$. In such a case, $\Gamma(s,w)$ is said to be a \emph{threshold representation of $G$}. For a threshold representation $\Gamma(s,w)$ of $G$, $s$ is called its threshold value, and $w$ is called its weight function. Note that we will always implicitly assume that the vertex set of $\Gamma(s,w)$ is the same as the domain of $w$.
\end{definition}
Let $\ell\geq 0$. If for graphs $G$ and $\{G_i\in \operatorname{TH}\}_{i\in [\ell]}$, we have $G=\bigcap_{i=1}^{\ell}G_i$, then we say that $\{G_i\in \operatorname{TH}\}_{i\in [\ell]}$ is a TH-representation of $G$.

It is known that $\operatorname{TH}\subset \operatorname{INT}$. This leads us to the following observation. \begin{observation}\label{thint}
    For any graph $G$, $\dimt(G)\geq \boxi(G)$.
\end{observation}
\noindent It is known that $\operatorname{co-{TH}}=\operatorname{TH}$. This gives us the following observation.
\begin{observation}\label{coTH}
    For any graph $G$, $\covt(G)=\dimt(\overline{G})$.
\end{observation}
\section{Kneser graphs as the intersection of $C_4$-free graphs}\label{kneser}
Throughout this section, we assume that $n\geq 2k+1$ and $k\geq 2$. Let 2K2F denote the class of $2K_2$-free graphs\footnote{$2K_2$ is the disjoint union of 2 edges.}. It is clear that 2K2F=co-C4F. In this section, we will prove Theorem~\ref{c4f} by showing that $\covk(\overline{K(n,k)})=n-2k$. 
\begin{lemma}\label{2K_2-freeupper}
    For $n\ge 2k+1$, $\covk(\overline{K(n,k)})\le n-2k$.
\end{lemma}
\begin{proof}
    We show a procedure to generate graphs $\{G_i\}_{i\in [n-2k]}$ such that $\overline{K(n,k)}=\bigcup_{i\in [n-2k]}G_i$. For any edge $\{A,B\}$ in $\overline{K(n,k)}$,
    \begin{itemize}
        \item if $(A\cup B)\cap [n-2k-1]$ is non-empty, assign it to $G_j$, where $j=\min ((A\cup B)\cap [n-2k-1])$;
        \item if $(A\cup B)\cap [n-2k-1]$ is empty, assign it to $G_{n-2k}$. 
    \end{itemize}This procedure ensures that $\overline{K(n,k)}=\bigcup_{i\in [n-2k]}G_i$. We will show that $\forall i\in [n-2k]$, $G_i$ is $2K_2$-free.
For $i\in [n-2k-1]$, let $\{A_1,B_1\}$ and $\{A_2,B_2\}$ be distinct edges in $G_i$. If $A_1,B_1,A_2$ and $B_2$ are not 4 distinct sets, then they cannot induce a $2K_2$ in $G_i$. So assume they are distinct. Let $X_1\in \{A_1,B_1\}$ and $X_2\in \{A_2,B_2\}$ be two sets that contain $i$. Then $X_1\cap X_2\neq \emptyset$ because it contains $i$. Moreover, $i=\min((X_1\cup X_2)\cap [n-2k-1])$. Therefore, there is an edge between $X_1$ and $X_2$ in $G_i$, and consequently $\{A_1,B_1,A_2,B_2\}$ does not induce a $2K_2$.

Now suppose $\{A_1,B_1\}$ and $\{A_2,B_2\}$ are distinct edges in $G_{n-2k}$. If $A_1,B_1,A_2$ and $B_2$ are not 4 distinct sets, then they cannot induce a $2K_2$ in $G_{n-2k}$. So assume they are distinct. The sets $A_1\cup B_1$ and $A_2\cup B_2$ cannot be disjoint because if they are, we would have  $|(A_1\cup B_1)\cup (A_2\cup B_2)|\geq 2k+2$, contradicting the fact that there is no element from $[n-2k-1]$ in $(A_1\cup B_1)\cup (A_2\cup B_2)$. Therefore, there are sets $X_1\in \{A_1,B_1\}$ and $X_2\in \{A_2,B_2\}$ such that $X_1\cap X_2\neq \emptyset$. Moreover, $X_1\cup X_2\subseteq[n]\setminus [n-2k-1]$. Therefore, $X_1$ and $X_2$ are adjacent in $G_{n-2k}$. 
\end{proof}
For a positive integer $d$, we will denote by $S^{d-1}$, the standard unit sphere in $\mathbb{R}^d$. We will use the following generalization of a theorem by Lusternik, Schnirelman and Borsuk, which was proved in \cite{Greene01122002}.
\begin{theorem}[\cite{Greene01122002}]\label{antipodal}
    If $S^{d-1}=C_0\cup C_1\cup  C_2\cup \cdots  \cup C_{d-1}$, such that for every $0\leq i\leq d-1$, $C_i$ is open or closed, then there is a $0\leq j\leq d-1$ such that $C_j$ contains two antipodal points.
\end{theorem}

\begin{lemma}\label{equator}
    There is an arrangement of any finite collection of points on the sphere $S^{d-1}$ so that no set of $d$ many points lies on an equator.
\end{lemma}
\begin{proof}
    Take the unit vectors in the direction defined by $(1,x,x^2, \ldots , x^{d-1})$. Any such point lying on a hyperplane, can be defined by a root of some polynomial of the form $c_0+c_1x+c_2x^2 +\cdots + c_{d-1}x^{d-1}$. Since a polynomial of degree $d-1$ can have at most $d-1$ roots, there can be at most $d-1$ such points on the hyperplane. Therefore, no set of $d$ points generated by this curve lies on an equator (an equator of the sphere is the intersection of the sphere with a hyperplane that passes through the origin). 
\end{proof}
\begin{lemma}\label{2K_2-free}
     For $n\ge 2k+1$, $\covk(\overline{K(n,k)})\ge n-2k$. 
\end{lemma}
\begin{proof}
      For $n= 2k+1$, the bound is trivially true. Therefore, assume that $n\geq 2k+2$. We show that any assignment of the edges of $\overline{K(n,k)}$ to at most $n-2k-1$ many $2K_2$-free graphs leads to a contradiction. Let $G_1, G_2, \ldots ,G_{n-2k-1}$ be a covering of the edges of $\overline{K(n,k)}$ such that for $i\in [n-2k-1]$, $G_i$ is a $2K_2$-free graph. 

      By Lemma~\ref{equator}, there is a way to embed $[n]$ into $S^{n-2k-1}$ such that no $n-2k$ points lie on an equator. Let $f$ be such an embedding. Any edge of $\overline{K(n,k)}$ is defined by two sets which have a common intersection. Therefore, distinct $k$-sets $A,B\in  \binom{[n]}{k}$ have an edge between them if and only if $k+1\leq |A\cup B|\leq 2k-1$.  
      
      For each $v\in S^{n-2k-1}$, define $H_v=\{x\in S^{n-2k-1}\colon \langle x,v \rangle> 0 \}$, which is an open hemisphere. Define the sets $C_0, C_1, C_2, \ldots ,C_{n-2k-1}$ as follows: for $v\in S^{n-2k-1}$, if there is a collection of points $T\subset f([n])\cap H_v$, and sets $A$ and $B$ that satisfy $k+1\leq |T|\leq 2k-1$ and $T=f(A)\cup f(B)$, then for every $i$ for which $G_i$ contains $\{A,B\}$ as an edge, put $v\in C_i$. Note that a point $v\in S^{n-2k-1}$ may be put into several sets $C_i$. Define $C_0=S^{n-2k-1}\setminus \bigcup_{i\in [n-2k-1]}C_i$. It is easy to see that for $i\in [n-2k-1]$, $C_i$ is open and that $C_0$ is closed (it is the complement of a union of open sets).  Therefore, we have $S^{n-2k-1}=C_0\cup \bigcup_{i\in [n-2k-1]}C_i$. By Theorem~\ref{antipodal}, there is a $0\leq j\leq n-2k-1$ such that $C_j$ contains two antipodal points $x_1$ and $x_2$. 
      
      Suppose $j\in [n-2k-1]$. Then by definition of $C_j$, there exist sets $T_1\subset f([n])\cap H_{x_1},T_2\subset f([n])\cap H_{x_2}$ and edges $\{A_1,B_1\}$ and $\{A_2,B_2\}$ in $G_j$ such that $k+1\leq |T_1|, |T_2|\leq 2k-1$ and $T_1=f(A_1)\cup f(B_1)$ and $T_2=f(A_2)\cup f(B_2)$. Since, $x_1$ and $x_2$ are antipodal points, $ H_{x_1}\cap  H_{x_2}=\emptyset$. Therefore, $(A_1\cup B_1)\cap (A_2\cup B_2)=\emptyset$. It follows that $\{A_1,B_2\},\{A_1,A_2\},\{B_1,B_2\},$ and $ \{B_1,A_2\}$ do not exist as edges in $\overline{K(n,k)}$. Therefore, the edges $\{A_1,B_1\}$ and $\{A_2,B_2\}$ induce a $2K_2$ in $G_j$. This is a contradiction. Therefore, we must have $j=0$. But this means that $|H_{x_1}\cap f([n])|\leq k$ and $|H_{x_2}\cap f([n])|\leq k$ (because for any $T$ subset of $k+1$ points of $f([n])$, we can find sets $A,B\in \binom{[n]}{k}$ such that $T=f(A)\cup f(B)$). This means that the equator defined by the plane $\langle x,x_1\rangle=0 $ contains at least $n-2k$ points. This contradicts the assumption on our embedding $f$. 
\end{proof}
\noindent The proof of Theorem~\ref{c4f} follows from Lemma~\ref{2K_2-freeupper} and Lemma~\ref{2K_2-free} and the fact that $\dimc(K(n,k))=\covk(\overline{K(n,k)})$.
\thbox*
\begin{proof}
    By Lemma~\ref{2K_2-free} the number of $2K_2$-free graphs needed to cover the edges of $\overline{K(n,k)}$ is at least $n-2k$. As a consequence, any covering of $\overline{K(n,k)}$ with at most $n-2k-1$ many co-interval graphs cannot be possible since co-interval graphs are also $2K_2$-free graphs. Hence, it follows that $\boxi(K(n,k))=\dimi(K(n,k))=\operatorname{cov_{co-INT}}(\overline{K(n,k)})\geq n-2k$ as claimed. 
\end{proof}

\section{Sublinear upper bounds on $\boxi(K(n,k))$ when $n\approx2k+1$}
Recall that $B(a,b;c)$ is the poset of height 2 consisting of the subsets of $[c]$ of size $a$ and $b$
ordered by the set-inclusion relation. 
 
\begin{lemma}\label{dimt}
   For $n\geq 2k+1$, $\dimt(K(n,k))\leq \dim (B(k-1,n-k-1;n-1))+\dimt(K(n-1,k))$.
\end{lemma}
\begin{proof}
We will consider $B_0$ to be the subposet of $B(k,n-k;n)$ obtained by only retaining sets that contain $1$. It is clear that $B_0$ is order isomorphic to $B(k-1,n-k-1;n-1)$. Let $\dim (B_0)=t_1$ and $\calL=\{L_i\}_{i\in [t_1]}$ be a realizer of size $t_1$ for $B_0$. For $i\in [t_1]$ and $X\in B_0 $ let $p_i(X)$ denote the position of $X$ in $L_i$. 

For $i\in [t_1]$, let $s_i=2\binom{n}{k}$ and define the weight functions $w_i\colon \binom{[n]}{k}\rightarrow \mathbb{Z}_{\geq 0}$ as follows: 

   $$w_i(Y)=\begin{cases}
        \binom{n}{k}+p_i([n]\setminus Y)\text{, if $Y$ does not contain 1}\\
       \binom{n}{k}-p_i(Y)\text{, if $Y$ contains 1}
    \end{cases}.$$ 

Let $G_0$ be the induced subgraph of $K(n,k)$ obtained by removing all sets that contain 1. It is clear that $G_0\cong K(n-1,k)$. Let $\dimt(G_0)=t_2$ and $\{\Gamma(s_j,w_j')\}_{j\in [t_2]}$ be the TH-representation of $G_0$. For every $j\in [t_2]$ and every $Y\in \binom{[n]}{k}$ let $w_j(Y)=\begin{cases}
    w_j'(Y), \text{ if $Y$ does not contain 1}\\
    s_j, \text{ if $Y$ does contain 1}
\end{cases}$.

We show that $\{\Gamma(s_i,w_i)\}_{i\in [t_1+t_2]}$ is a TH-representation of $K(n,k)$. It is clear that the above graphs are super graphs of $K(n,k)$. We consider the following three cases. 
\begin{itemize}
\item If $Y_1,Y_2\in \binom{[n]}{k}$ be such that both of them contain $1$, then they are not adjacent. Clearly for any $i\in [t_1]$, $w_i(Y_1)+w_i(Y_2)=\binom{n}{k}-p_i(Y_1)+\binom{n}{k}-p_i(Y_2)<2\binom{n}{k}$. 
\item If $Y_1,Y_2\in \binom{[n]}{k}$ be such that both of them do not contain $1$, then they are adjacent in $K(n,k)$ if and only if they are adjacent in $G_0$.  Clearly $Y_1,Y_2$ are adjacent in $G_0$ if and only if $\forall j\in [t_2]$, $w_j(Y_1)+w_j(Y_2)=w_j'(Y_1)+w_j'(Y_2)\geq s_j$. 
\item If $Y_1,Y_2\in \binom{[n]}{k}$ be such $Y_1$ contains 1 and $Y_2$ does not contain $1$, then they are adjacent if and only if $Y_1\subset [n]\setminus Y_2$, which holds if and only if $\forall i\in [t_1]$, $p_i(Y_1)< p_i([n]\setminus Y_2)$ which in turn holds if and only if  $\binom{n}{k}-p_i(Y_1)+ \binom{n}{k}+p_i([n]\setminus Y_2)\geq 2\binom{n}{k}$. Therefore, $Y_1$ and $Y_2$ are adjacent if and only if $\forall i\in [t_1]$ $w_i(Y_1)+w_i(Y_2)\geq s_i$.
\end{itemize}
\end{proof}
Note that the bound in Lemma~\ref{dimt} is tight up to a constant factor as $\dimt(K(n,k))\ge \dimt(K(n-1,k))$ and $\dimt(K(n,k))\ge \boxi(K(n,k))\ge \frac{1}{2}\dim (B(k-1,n-k-1;n-1))-2$. 

\noindent Now we are ready for the proof of Theorem~\ref{sublinear}.
\sublinear*
\begin{proof}
    Let $n=2k+t$ in the statement of Lemma~\ref{dimt}. By Corollary 1.12 of~\cite{BKKT}, $\dim( B(k-1,k+t-1;2k+t-1))\leq O(t^2\log (2k+t-1))$. Solving the resulting recursion gives $\dimt(K(2k+t,k))\leq  O( t^3\log(2k+t) )+\dimt(K(2k,k))\leq  O( t^3\log(2k+t))$ where the last inequality follows by observing that $\dimt(K(2k,k))\leq 2$ (the graph $K(2k,k)$ is a disjoint union of edges). The statement then follows by Observation~\ref{thint}.
\end{proof}
\odd*
\begin{proof}
    The proof of this follows exactly as the proof of Theorem~\ref{sublinear}, except that we use Theorem 1 of~\cite{kostochka1997dimension} instead of Corollary 1.12 of~\cite{BKKT}. 
\end{proof}
\textbf{Declaration of AI use:} AI tools were not used to write this paper, and all the proofs are due to the authors. 

\textbf{Acknowledgement:} The authors thank Rishikesh Gajjala for the fruitful discussions during the initial stages of this work.
\bibliographystyle{plain}
\bibliography{ref}
\end{document}